\documentclass[review,onefignum,onetabnum]{siamart190516}

\usepackage{graphicx}
\usepackage{cite}
\usepackage{caption}
\usepackage{subcaption}
\usepackage{epstopdf}
\usepackage{amsmath}
\usepackage{amsfonts}
\usepackage{cleveref}

\newtheorem{remark}{Remark}
\nolinenumbers
\headers{Stabilization of production model}{S. G\"{o}ttlich, M. Herty, and G. Weldegiyorgis}

\title{Input-to-state stability of a scalar conservation law with nonlocal velocity}

\author{Simone G\"{o}ttlich\thanks{Department of Mathematics, University of Mannheim, Germany (\email{goettlich@unimannheim.de}).}
	\and Michael Herty\thanks{IGPM, Templergraben 55, RWTH Aachen University, Germany (\email{herty@igpm.rwth-aachen.de}).}
	\and Gediyon Weldegiyorgis\thanks{Department of Mathematics and Applied Mathematics, University of Pretoria, South Africa (\email{gediyon@aims.ac.za}).}}

\begin{document}	
\maketitle 

\begin{abstract}
	In this paper, we study input-to-state stability (ISS) of an equilibrium for a scalar conservation law with nonlocal velocity and measurement error arising in a highly re-entrant manufacturing system. By using a suitable Lyapunov function, we prove sufficient and necessary conditions on ISS. We also analyze the numerical discretization of ISS for a discrete scalar conservation law with nonlocal velocity and measurement error. A suitable discretized Lyapunov function is also analyzed to provide ISS of an equilibrium for the numerical approximation. Finally, we show numerical simulations to validate the theoretical findings. 
\end{abstract}



\section{Introduction}\label{Intro}
\addcontentsline{toc}{section}{Introduction} 

The nature of modern high-volume production is characterized by a large number of items passing through many production steps. This type of production system has fluid-like properties and has been modelled successfully by continuum models \cite{Armbruster2006, Herty2007, Armbruster2006a, DApice2010,Chen1988}. In these models, the product at different production stages and the speed of production are the quantities of interest. 

Specifically, in the manufacturing system of a factory that involves a highly re-entrant system where products visit machines multiple times, such as the production of semiconductor devices, a continuum model has been introduced in \cite{Armbruster2006a} that is inspired by the Lighthill – Whitham traffic model \cite{Helbing1996}. The dynamics of this model is mathematically given by hyperbolic partial differential equation of the form 
\begin{equation}\label{eq:Continuum_model}
	\partial_t\rho(t,x) + \lambda(W(t))\partial_x\rho(t,x) = 0,\quad t \in [0, +\infty),\; x \in [0,1], 
\end{equation}
where $ \rho(t,x) $ is the product density which describes the total mass $W(t)$ at the time $ t $ and the production stage $ x $, 
\begin{equation}\label{eq:Total_mass}
	W(t) = \int_{0}^{1} \rho(t,x)dx,\quad  t \in (0, +\infty). 
\end{equation}
 Contrary to classical traffic flow models the differential equation depends on the {\bf nonlocal} quantity \eqref{eq:Total_mass}. The function  $\lambda(W(t)) $ is a velocity. In production systems, it is natural to assume that the  velocity function is positive and decreasing as the total mass is increasing. In the manufacturing system, the initial density of products at production stage $x$ is taken as the initial data 
\begin{equation}\label{eq:Continuum_model_IC}
\rho(0,x) = \rho_0(x),\quad x \in [0,1],
\end{equation} 
and the influx is used to control the system or stabilize the system at an equilibrium. Since the velocity is positive, we only require 
boundary conditions at $ x = 0 $ i.e. the influx
\begin{equation}\label{eq:Continuum_model_BC1}
\rho(t,0)\lambda(W(t)) = U(t),\quad t \in [0, +\infty). 
\end{equation}
Under suitable assumptions on $\lambda, \rho_0$ and $U$,  the existence and uniqueness of a classical solution of the Cauchy problem for the scalar conservation law \cref{eq:Continuum_model} with \cref{eq:Continuum_model_IC,eq:Continuum_model_BC1}  is proven in \cite{Coron_2013, Coron2010, Shang2011, Chen_2017}. 

General stabilization problems with boundary controls have been studied in the past years in \cite{Bastin2008,Tanwani_2018,Bastin,Bastin2011,Bastin2017,Bastin_2016,Coron2004,Diagne2012,Gugat2018,LiHui2009,Prieur2012,Coron2008, Zhang2017} for hyperbolic systems and recently in \cite{Chen_2017,Coron_2013} for scalar conservation laws with nonlocal velocity. The focus is to derive an asymptotic stability around a given equilibrium such that solutions to the conservation laws  reach the equilibrium state as time tends to infinity. Such a property is attained by an exponential stability result e.g. in \cite{Bastin_2016}. However, when boundary controls are subjected to unknown disturbances, solutions reaching the given equilibrium point are influenced by the disturbances and a notion of asymptotic stability is required.  This property is covered in an input-to-state stability (ISS) \cite{Lamare2018,Prieur2012,Tanwani_2018}. Concerning an asymptotic behavior of classical solutions, the Lyapunov method is used to investigate sufficient conditions to achieve an exponential stability in \cite{Bastin2008,Coron2004,Diagne2012} for hyperbolic systems and in \cite{Chen_2017,Coron_2013} for scalar conservation laws with nonlocal velocity. The Lyapunov  method is also used for ISS of (local) hyperbolic systems in \cite{Tanwani_2018,Prieur2012}. For the numerical analysis of asymptotic behavior of numerical solutions discretized by a first-order finite volume scheme, a discrete Lyapunov function is used to prove exponential stability results for hyperbolic systems in \cite{Banda,Banda2013,Goettlich2017,Gerster2019,Goettlich2017a} and for scalar conservation laws with nonlocal velocity in \cite{Chen_2017}, and ISS results for (local) hyperbolic systems could be established recently in \cite{Weldegiyorgis2019a, Weldegiyorgis, Bastin2020}.       

In connection with a scalar conservation law with nonlocal velocity, in \cite{Chen_2017}, the authors have studied global feedback stabilization of the closed-loop system \cref{eq:Continuum_model} under the feedback law
\begin{equation}\label{eq:feedback_law1}
U(t) - \rho^*\lambda(\rho^*) = k\Big(\rho(t,1)\lambda(W(t)) - \rho^*\lambda(\rho^*)\Big),\; t \in (0, +\infty),
\end{equation}
where $ k\in \mathbb{R} $ is the feedback parameter and $ \rho^* \in \mathbb{R} $ is a given equilibrium. They generalize the stabilization results of \cite{Coron_2013} by using a Lyapunov function. In particular, for a given equilibrium $ \rho^* = 0 $ and a general velocity function $ \lambda \in C^1([0, +\infty);[0, +\infty))$, the global stabilization result in $ L^2 $ for the closed-loop system \cref{eq:Continuum_model,eq:Continuum_model_IC,eq:feedback_law1} is generalized to $ L^p $ $ (p \geq 1) $. Then, the global stabilization result in $ L^2 $ for the closed-loop system \cref{eq:Continuum_model,eq:Continuum_model_IC,eq:feedback_law1} with a family of velocity functions 
\begin{equation}\label{eq:Velocity_function}
\lambda(s) = \frac{A}{B + s}, \quad s \in [0, +\infty) \quad \text{with}\quad A > 0\; B>0, 
\end{equation}  
is obtained for a given equilibrium $ \rho^* > 0 $. By using a discrete Lyapunov function, they also established stabilization results for a discrete scalar conservation law with nonlocal velocity and using  a first-order finite volume scheme.  

In this paper, we study ISS for the closed-loop system  \cref{eq:Continuum_model,eq:Continuum_model_IC} under the feedback law defined by
\begin{equation}\label{eq:feedback_law2}
U(t) - \rho^*\lambda(\rho^*) = k\Big(\left(\rho(t,1) + d(t)\right)\lambda(W(t)) - \rho^*\lambda(\rho^*)\Big),\; t \in (0, +\infty),
\end{equation}   
where $ d \in \mathbb{R}$ is a bounded perturbation in the measurement. In particular,  we use an ISS-Lyapunov function to investigate sufficient and necessary conditions for ISS in $ L^2 $ for an equilibrium $ \rho^* \geq 0 $ and the velocity function defined by \cref{eq:Velocity_function}.  The numerical analysis of sufficient and necessary conditions for ISS is performed by using a discrete ISS-Lyapunov function for numerical solution obtained by a first--order finite volume scheme. Moreover, we provide numerical simulations to illustrate theoretical results for some velocity functions of type \cref{eq:Velocity_function}.   

The paper is organized as follows: In \cref{sec:Section01}, we present stabilization results of ISS for a scalar conservation law with nonlocal velocity and measurement error. The numerical discretization of stabilization results of ISS for the scalar conservation law with nonlocal velocity and measurement error is presented in \cref{sec:Section02}. Finally, in \cref{sec:Section03}, we show numerical simulations for the scalar conservation law with nonlocal velocity and measurement error to illustrate the theoretical results.     

\section{Asymptotic stability of a scalar conservation law with nonlocal velocity and measurement error}\label{sec:Section01}

We study ISS of a closed-loop system of scalar conservation laws with nonlocal velocity and measurement error of the form:
\begin{equation}\label{eq:Nonlocal_system}
	\begin{cases}
	\partial_t\rho(t,x) + \lambda(W(t))\partial_x\rho(t,x) = 0,\quad t \in (0, +\infty),\; x \in (0,1),\\
	\rho(0,x) = \rho_0(x),\quad x \in (0,1),\\
	U(t) - \rho^*\lambda(\rho^*) = k(\left(\rho(t,1) + d(t)\right)\lambda(W(t)) - \rho^*\lambda(\rho^*)),\quad t \in (0, +\infty), \\
	W(t) = \int_{0}^{1} \rho(t,x)dx,\quad t \in (0, +\infty),
	\end{cases}
\end{equation}
where $ \rho(t,x) $ is the product density, $ \lambda(\cdot) \in C^1([0,+\infty), (0,+\infty))$ is the velocity function, $ W(t) $ is total mass defined by \cref{eq:Total_mass}, $ U(t) $ is the controller or the input function defined by \cref{eq:Continuum_model_BC1}, $ k \in [0,1] $ is a feedback parameter, $ \rho^* \geq 0 $ is an equilibrium solution and $ d(t) \in \mathbb{R}$ is a bounded (known) perturbation in the measurement.  A weak solution of the closed-loop system \cref{eq:Nonlocal_system} is defined below.
\begin{definition}[Weak solution]\label{def:Nonlocal_weak_solution}
	A function $ \rho \in C^0([0,T];L^1(0,1)) $ is called a weak solution to \cref{eq:Nonlocal_system} if for every $ T > 0 $, every $ s \in (0, T] $ and every $ \varphi \in C^1([0,s]\times[0,1])$ satisfying 
	\begin{equation*}
	\varphi(s, x) = 0,\; \forall x \in [0,1]\quad \text{and}\quad \varphi(t, 1) = \kappa\varphi(t, 0),\; \forall t \in [0,s], 
	\end{equation*}
	the following equation holds: 
	\begin{align*}
	&\int_{0}^{s}\int_{0}^{1} \rho(t,x)\left(\partial_t\varphi(t,x) + \lambda(W(t))\partial_x\varphi(t,x)\right)dxdt \\
	& + \int_{0}^{s} ((1- k) \rho^*\lambda(\rho^*) + d(t))\varphi(t, 0)dt + \int_{0}^{1} \rho(0,x)\varphi(0,x)dx = 0.
	\end{align*}
\end{definition}

Let $d\equiv 0$, $ \rho^* \geq 0 $, $ p \in [1,+\infty) $, and $ k\in [0,1] $ be given. Then, the existence and uniqueness of the non-negative weak solution $ \rho \in C^0([0, +\infty);L^p(0,1)) $ and the non-negative classical solution $ \rho \in C^1([0, +\infty)\times[0,1]) $ of the closed-loop system \cref{eq:Nonlocal_system} are available in \cite{Chen_2017,Coron_2013}.  

We now analyze ISS for the system \cref{eq:Nonlocal_system} with $ \rho^* \geq 0 $ in the sense of the following definitions:

\begin{definition}[Input-to-state stability (ISS)]\label{def:ISS} 
	An equilibrium $\rho^* \geq 0 $ of the closed-loop system \cref{eq:Nonlocal_system} is an ISS in $L^2-$norm with respect to a bounded disturbance function $t \to d(t) $ if there exist positive constants $\gamma_1 > 0$, $ \gamma_2 > 0 $ and $ \gamma_3 > 0 $ such that, for every initial condition $\rho_0(x) \in L^2$, the  $L^2-$solution to the closed-loop system \cref{eq:Nonlocal_system} satisfies
	\begin{equation}\label{eq:ISS_condition}
	{\|\rho(t,\cdot) \|}_{L^2} \leq \gamma_2 {e}^{-\gamma_1 t}{\|\rho_0 \|}_{L^2} + \gamma_3 \sup_{s \in [0,t]}\left(|d(s)|\right),\;t \in [0, +\infty).
	\end{equation}	
\end{definition}

\begin{definition}[ISS-Lyapunov function]\label{def:ISS_Lyapunov_function}
	A continuously differentiable function $ \mathcal{L}:[0,\infty) \rightarrow \mathbb{R}_+ $ is said to be an ISS-Lyapunov function for the closed-loop system \cref{eq:Nonlocal_system} if 
	\begin{itemize}
		\item[(i)] there exist positive constants $ \alpha_1 >0$ and $\alpha_2 >0 $ such that for all solutions $ \rho $ and $ t \in [0, +\infty)$
		\begin{equation}\label{eq:ISS_Lyapunov_function_condition1}
		\alpha_1\|\rho(t,\cdot) \|_{L^2}^2 \leq \mathcal{L}(t) \leq \alpha_2\|\rho(t,\cdot) \|_{L^2}^2,
		\end{equation}
		\item[(ii)] there exist  positive constants $ \eta > 0$ and $ \nu > 0 $ such that for all solutions $ \rho $ and $ t \in [0, +\infty)$ 
		\begin{equation}\label{eq:ISS_Lyapunov_function_condition2}
		\mathcal{L}'(t) \leq -\eta \mathcal{L}(t) + \nu d^2(t).
		\end{equation}
	\end{itemize}
\end{definition}

\begin{theorem}[ISS for $ \rho^* \geq 0 $]\label{thm:ISS_nonlocal2}
	For every $\rho^* \geq 0$, every $ k \in [0,1) $, every $ R > 0 $ and for every $ \rho_0 \in L^2(0,1) $ satisfying $ \rho_0 \geq 0 $ a.e. in $ (0,1) $ and 
	\begin{equation}\label{eq:ISS_nonlocal2_condition1}
	\|\rho_0(\cdot) - \rho^*\|_{L^2(0,1)} \leq R.
	\end{equation} 
	Assume there exists a unique, non-negative almost everywhere weak solution \\
	 $ \rho \in C^0([0, +\infty); L^2(0,1))$ to the Cauchy problem \cref{eq:Nonlocal_system} with \cref{eq:Velocity_function}. Then, the steady-state $ \rho^* $ of the system \cref{eq:Nonlocal_system} is ISS in $L^2-$norm with respect to the disturbance function $ d $.
\end{theorem}
\begin{remark}
	For $ \rho^* = 0 $, any velocity function $ \lambda(\cdot) \in C^1([0,+\infty), (0,+\infty))$ can be considered in \cref{thm:ISS_nonlocal2} (see Theorem 3.1 in \cite{Chen_2017} for details of exponential stability of the system \cref{eq:Nonlocal_system} when $ d \equiv  0 $). 
\end{remark}

Before we begin the proof of \cref{thm:ISS_nonlocal2}, we consider the following transformation around the equilibrium $ \rho^* $, 
\begin{equation*}
\begin{split}
&\tilde{\rho}(t,x) := \rho(t,x) - \rho^*, \quad \widetilde{W}(t) := W(t) - \rho^*, \quad \tilde{\rho_0}(x) := \rho_0(x) - \rho^*, \\
&\tilde{\lambda}_{\widetilde{W}}(t) := \lambda(\rho^* + \widetilde{W}(t)),\quad \widetilde{U}(t) := \tilde{\lambda}_{\widetilde{W}}(t) \tilde{\rho}(t,0).
\end{split}
\end{equation*}
Then, the system \cref{eq:Nonlocal_system} with \cref{eq:Velocity_function} can be rewritten as follows for $t \in (0, +\infty)$: 
\begin{equation}\label{eq:Nonlocal_system2}
\begin{cases}
\partial_t\tilde{\rho}(t,x) + \tilde{\lambda}_{\widetilde{W}}(t)\partial_x\tilde{\rho}(t,x) = 0, \; x \in (0,1),\\
\tilde{\rho}(0,x) = \tilde{\rho}_0(x),\; x \in (0,1), \\
\widetilde{U}(t) = k\tilde{\lambda}_{\widetilde{W}}(t) \left(\tilde{\rho}(t,1) + d(t)\right) + (1-k)\rho^*\left(\lambda(\rho^*) -  \tilde{\lambda}_{\widetilde{W}}(t)\right), \\
\tilde{\lambda}_{\widetilde{W}}(t) := \lambda(\rho^* + \widetilde{W}(t)), \\
\widetilde{W}(t) = \int_{0}^{1} \tilde{\rho}(t,x)dx \geq - \rho^*, \\
\lambda(s) = \frac{A}{B + s},\quad \text{with} \quad  A > 0,\; B > 0,\; s \in [0, +\infty).
\end{cases}
\end{equation}
By using the velocity function \cref{eq:Velocity_function} in \cref{eq:Nonlocal_system2}, we have   
\begin{equation}\label{eq:Boundary_term0} 
\rho^*\left(\lambda(\rho^*) - \tilde{\lambda}_{\widetilde{W}}(t) \right) = {\theta}\tilde{\lambda}_{\widetilde{W}}(t)\widetilde{W}(t), \quad t \in [0, +\infty), 
\end{equation}
where $\displaystyle \theta:= \frac{\rho^*}{B + \rho^*}$. 

For convenience, until the end of this proof, we omit the symbol ``\textasciitilde". Then, the system \cref{eq:Nonlocal_system2} with \cref{eq:Boundary_term0}  can be rewritten in the following form for $t \in (0, +\infty)$:
\begin{equation}\label{eq:Nonlocal_system3}
\begin{cases}
\partial_t\rho(t,x) + \lambda_{W}(t)\partial_x\rho(t,x) = 0,\; x \in (0,1),\\
\rho(0,x) = \rho_0(x),\; x \in (0,1), \\
U(t) = k{\lambda}_{W}(t) \left(\rho(t,1) + d(t)\right) + (1 - k){\theta}\lambda_{W}(t)W(t)\;\text{with}\; \theta = \frac{\rho^*}{B + \rho^*}, \\
\lambda_{W}(t) := \lambda(\rho^* + W(t)), \\
W(t) = \int_{0}^{1} \rho(t,x)dx \geq - \rho^*, \\
\lambda(s) = \frac{A}{B + s},\quad \text{with} \quad  A > 0,\; B > 0,\; s \in [0, +\infty).
\end{cases}
\end{equation}
With the above notation, the assumption \cref{eq:ISS_nonlocal2_condition1} in \cref{thm:ISS_nonlocal2} is then
\begin{equation}\label{eq:ISS_nonlocal2_condition1_updated}
\|\rho_0\|_{L^2(0,1)} \leq R.
\end{equation}
\begin{proof} 
	The following proof of \cref{thm:ISS_nonlocal2} is an extension of the proof of Theorem 3.2 in \cite{Chen_2017}. Since $ C^1-$functions are dense in $ L^2(0,1) $, we can analyze ISS for the system \cref{eq:Nonlocal_system3} with non-negative weak solution $ \rho \in C^0([0, +\infty);L^2(0,1))$ as follows:
	We first define a candidate ISS-Lyapunov function by 
	\begin{equation}\label{eq:Lyapunov_function}
	\mathcal{L}(t) := \int_{0}^{1} \rho^2(t,x)e^{-\beta x}dx + a W^2(t), \quad \forall t \in [0, +\infty), 
	\end{equation} 
	where $ \beta > 0 $ and $ a \in \mathbb{R} $ are constants. If
	\begin{equation}\label{eq:Constant_a} 
	a > -\frac{\beta}{e^{\beta}-1},
	\end{equation} 
	then $ \mathcal{L}(t) > 0 $ for all $ t \geq 0 $ and 
	there exist positive constants $ C_i = C_i(a, \beta)$, $i \in \{1,2\} $ such that for all $ t \geq 0 $
	\begin{equation}\label{eq:Lyapunov_approximation1} 
	W^2(t) \leq C_1 \int_{0}^{1} \rho^2(t,x)e^{-\beta x}dx \leq  \mathcal{L}(t) \leq C_2 \int_{0}^{1} \rho^2(t,x)e^{-\beta x}dx.  
	\end{equation}
	In particular, for $\rho^* = 0  $, we take $ a = 0 $ in \cref{eq:Lyapunov_function}. Then, the time derivative of the candidate ISS-Lyapunov function \cref{eq:Lyapunov_function} is computed as follows: 
	\small
	\begin{align}
	\frac{d\mathcal{L}}{dt} =&\; \int_{0}^{1} 2\rho(t,x)\rho_t(t,x)e^{-\beta x}dx + 2a W(t)\frac{dW}{dt},\notag\\
	=&\; -\beta\lambda_{W}(t)\int_{0}^{1} \rho^2(t,x)e^{-\beta x}dx\notag\\
	& + \frac{1}{\lambda_{W}(t)} \left(U^2(t) - \left(\lambda_{W}(t)\rho(t,1)\right)^2e^{-\beta}\right) + 2a W(t) \left(U(t) - \lambda_{W}(t)\rho(t,1) \right), \notag\\
	=&\; -\beta\lambda_{W}(t)\int_{0}^{1} \rho^2(t,x)e^{-\beta x}dx \notag\\
	& + \frac{1}{\lambda_{W}(t)} \left(\left[k\lambda_{W}(t)\rho(t,1) + (1 - k){\theta}\lambda_{W}(t)W(t) + k\lambda_{W}(t)d(t)\right]^2 - \left(\lambda_{W}(t)\rho(t,1)\right)^2e^{-\beta}\right) \notag\\
	& + 2a W(t) \left(\left[k\lambda_{W}(t)\rho(t,1) + (1 - k){\theta}\lambda_{W}(t)W(t) + k\lambda_{W}(t)d(t)\right] - k\lambda_{W}(t)\rho(t,1) \right),\notag\\
	\leq &\; -\beta\lambda_{W}(t)\int_{0}^{1} \rho^2(t,x)e^{-\beta x}dx  + 3k^2\lambda_{W}(t)d^2(t) + \lambda_{W}(t)b(t),\label{eq:Lyapunov_approximation2} 
	\end{align}
	\normalsize
	where the boundary term is defined as
	\begin{align*}
	b(t) &:= 2\left[k\rho(t,1) + (1 - k){\theta}W(t)\right]^2 - \rho^2(t,1)e^{-\beta}\notag\\
	&\quad + 2a(k-1)W(t) \left(\rho(t,1) - {\theta}W(t)\right) + a^2W^2(t).
	\end{align*} 
	The boundary term can be further simplified as
	 
	\begin{align*}
	b(t) =&\; \left(2k^2 -  e^{-\beta}\right)\rho^2(t,1) + 2(2k\theta - a) (1-k) \rho(t,1) W(t)\notag\\
	      &\quad + \left(2(1-k)^2\theta^2 + 2a(1-k)\theta + a^2 \right)W^2(t),\\
	=&\; \left(2k^2 -  e^{-\beta}\right)\left[\rho(t,1) + \frac{(2k\theta - a) (1-k)}{2k^2 -  e^{-\beta}} W(t)\right]^2 - \frac{(2k\theta - a)^2 (1-k)^2}{2k^2 -  e^{-\beta}}W^2(t)\\
	 &\quad \left(2(1-k)^2\theta^2 + 2a(1-k)\theta  + a^2 \right)W^2(t),\\
	= &\; \left(2k^2 -  e^{-\beta}\right)\left[\rho(t,1) + \frac{(2k\theta - a) (1-k)}{2k^2 -  e^{-\beta}} W(t)\right]^2\\
	&\quad - \frac{(1-k)^2}{2k^2 -  e^{-\beta}}\left[(2k\theta - a)^2 - 2\theta^2\left(2k^2 -  e^{-\beta}\right) - \frac{2a\theta\left(2k^2 -  e^{-\beta}\right)}{1-k} \right]W^2(t) 
	 + a^2 W^2(t),\\
	= &\; \left(2k^2 -  e^{-\beta}\right)\left[\rho(t,1) + \frac{(2k\theta - a) (1-k)}{2k^2 -  e^{-\beta}} W(t)\right]^2\\
	&\quad - \frac{(1-k)^2}{2k^2 -  e^{-\beta}}\left[a  - \theta\left(\frac{2k -  e^{-\beta}}{1-k}\right) \right]^2W^2(t) \\
	&\quad + \left[\frac{\theta^2\left(2k -  e^{-\beta}\right)^2}{2k^2 -  e^{-\beta}} - \frac{2\theta^2(1-k)^2e^{-\beta}}{2k^2 -  e^{-\beta}} + a^2\right] W^2(t),\\
	= &\; \left(2k^2 -  e^{-\beta}\right)\left[\rho(t,1) + \frac{(2k\theta - a) (1-k)}{2k^2 -  e^{-\beta}} W(t)\right]^2\\
	&\quad - \frac{(1-k)^2}{2k^2 -  e^{-\beta}}\left[a  - \theta\left(\frac{2k -  e^{-\beta}}{1-k}\right) \right]^2W^2(t)  +   \left[\theta^2\left(2 -  e^{-\beta}\right) + a^2\right] W^2(t). 
	\end{align*}
	\normalsize
	For given $ k \in [0,1) $ and $ \rho^* > 0 $, if we choose $ \beta > 0 $ such that 
	\begin{equation}\label{eq:Beta_condition1} 
	e^{-\beta} > 2k \geq 2k^2, 
	\end{equation} 
	and take 
	\begin{equation}\label{eq:Constant_a2} 
	a := \theta \left(\frac{2k - e^{-\beta}}{1 -k} \right) \leq 0,
	\end{equation} 	
	we obtain that 
	\begin{equation}\label{eq:Boundary_term_approx} 
	b(t) \leq \left[\left(2 -  e^{-\beta}\right) + \left(\frac{2k - e^{-\beta}}{1 -k} \right)^2\right]{\theta}^2W^2(t).
	\end{equation}  	 
	Therefore, from \cref{eq:Lyapunov_approximation1,eq:Lyapunov_approximation2,eq:Boundary_term_approx}, we get 
	\small
	\begin{align}
	\frac{d\mathcal{L}}{dt} &\leq -\left[\beta - {\theta}^2\left(\frac{e^{\beta}-1}{\beta}\right)\left(\left(2 -  e^{-\beta}\right) + \left(\frac{2k - e^{-\beta}}{1 -k} \right)^2\right)\right]\lambda_{W}(t)\int_{0}^{1} \rho^2(t,x)e^{-\beta x}dx\notag \\
	&\quad + \frac{3}{2}e^{-\beta}\lambda_{W}(t)d^2(t).\label{eq:Lyapunov_approximation3} 
	\end{align}	
	\normalsize
	By taking $ \beta = \beta(\rho^*,k) > 0 $ such that \cref{eq:Constant_a}  and \cref{eq:Beta_condition1}  hold and 
	\begin{equation}\label{eq:Beta_condition2} 
	\beta - {\theta}^2\left(\frac{e^{\beta}-1}{\beta}\right)\left(\left(2 -  e^{-\beta}\right) + \left(\frac{2k - e^{-\beta}}{1 -k} \right)^2\right) > 0,	
	\end{equation}
	and from \cref{eq:Lyapunov_approximation1}, we have 
	\begin{equation}\label{eq:Lyapunov_approximation4} 
	\frac{d\mathcal{L}}{dt} \leq  - f(\beta)\lambda_{W}(t)\mathcal{L}(t) + \frac{3}{2}e^{-\beta}\lambda_{W}(t)d^2(t),
	\end{equation}
	where
	\begin{equation*}
	f(\beta) = \frac{1}{C_1}\left[\beta - {\theta}^2\left(\frac{e^{\beta}-1}{\beta}\right)\left(\left(2 -  e^{-\beta}\right) + \left(\frac{2k - e^{-\beta}}{1 -k} \right)^2\right)\right]. 
	\end{equation*}
	We solve \cref{eq:Lyapunov_approximation4} to obtain  
	\begin{align}
	\mathcal{L}(t) &\leq e^{-\int_{0}^{t}f(\beta)\lambda_{W}(s)ds}\mathcal{L}(0) + \frac{3}{2}e^{-\beta}\int_{0}^{t}\lambda_{W}(s)d^2(s) e^{-\int_{s}^{t}f(\beta)\lambda_{W}(r)dr}ds,\notag\\
	&\leq e^{-\int_{0}^{t}f(\beta)\lambda_{W}(s)ds}\mathcal{L}(0) + \frac{3}{2}\frac{e^{-\beta}}{f(\beta)} d^2(t),\quad t \in [0, +\infty).\label{eq:Lyapunov_approximation5} 
	\end{align} 	
	Moreover, from \cref{eq:Nonlocal_system3,eq:ISS_nonlocal2_condition1_updated,eq:Lyapunov_approximation1,eq:Lyapunov_approximation5}, there exist positive constants $ C_i = C_i(\beta)$, $i \in \{3,4\} $ such that 
	\begin{equation}\label{eq:Lyapunov_approximation6} 
	- \rho^* \leq W(t) \leq \sqrt{C_2C_3R^2 + C_4\sup_{s \in [0,t]}(d^2(s))}, \quad \forall t \in [0, +\infty).
	\end{equation}
	Let 
	\begin{equation}\label{eq:Lyapunov_approximation7} 
	\sigma_2 := \inf_{W(t)} \lambda_{W}(t) > 0,\quad\text{and}\quad \delta_2 := \sup_{W(t)} \lambda_{W}(t).
	\end{equation}
	Therefore, from \cref{eq:Lyapunov_approximation1,eq:Lyapunov_approximation5,eq:Lyapunov_approximation7}, for all $ t \in [0, +\infty) $, we have 
	\begin{equation}\label{eq:Lyapunov_approximation8} 
	\begin{split}
	C_1\|\rho(t,\cdot) \|_{L^2(0,1)}^2 \leq \mathcal{L}(t) &\leq  e^{-\sigma_2 f(\beta) t}\mathcal{L}(0) + \frac{3}{2}\frac{\delta_2e^{-\beta}}{\sigma_2f(\beta)}\sup_{s \in [0,t]}\left(d^2(s)\right),\\
	&\leq C_2e^{-\sigma_2 f(\beta) t}\|\rho_0 \|_{L^2(0,1)}^2 + \frac{3}{2}\frac{\delta_2e^{-\beta}}{\sigma_2f(\beta)}\sup_{s \in [0,t]}\left(d^2(s)\right). 
	\end{split}
	\end{equation}
	Thus, we completed the proof of \cref{thm:ISS_nonlocal2}.	
\end{proof}

\section{Numerical study of asymptotic stability of a scalar conservation law with nonlocal velocity and measurement error}\label{sec:Section02}

In order to study ISS of the closed-loop system \cref{eq:Nonlocal_system}, we first divide the spatial domain $ [0,1] $ with the cells centers and boundary points are denoted by $ x_j = (j-\frac{1}{2}){\Delta x},\;j \in \{1,\ldots,J\} $ and, $ x_0 $ and $ x_{J+1} $, respectively such that $ {\Delta x}J = 1,\; J \in \mathbb{N} $. Moreover, we approximate $ W(t) $ as 
\begin{equation}\label{eq:Disc_total_mass}
W^n = {{\Delta x}}\sum_{j = 1}^{J}\rho_{j}^{n},\quad n \in \{1,2,\ldots\}, 
\end{equation}
with the point wise values of the solution $ \rho_{j}^{n} = \rho(t^n, x_j)$ and the discrete velocity function $ \lambda^{n} $ is given by
\begin{equation}\label{eq:Disc_velocity_function}
	\lambda^n = \lambda(W^n) = \frac{A}{B + W^n},\quad A > 0,\; B > 0,
\end{equation}
where $ t^n = n {\Delta t},\; n \in \{0, 1, \ldots \}$ denotes the discrete time such that the time step size $ {\Delta t} $ satisfies the CFL condition given by
\begin{equation}\label{eq:CFL_cond}
	0 < r^n:= \frac{\lambda^n {{\Delta t}}}{{\Delta x}} \leq 1,\; \forall n \in \{0,1,\ldots\}.
\end{equation}

For given initial values $ \vec{\rho}^0 = (\rho_{0}^{0},\rho_{1}^{0}, \dots, \rho_{J}^{0})^\top $ with $ \rho_{j}^{0} \geq 0$, $j \in \{0, \dots, J\} $, we use a first--order finite volume scheme, given by the explicit Upwind method, to discretize the system \cref{eq:Nonlocal_system} with $ \rho^* \geq 0 $ as follows:    
\begin{equation}\label{eq:Disc_nonlocal_system1}
\begin{cases}
\rho_{j}^{n+1} = \left( 1 - r^n\right) \rho_{j}^{n} + r^n \rho_{j-1}^{n}, & j \in \{1, \dots, J\},\; n \in \{0,1,\ldots\},\\
\rho_{0}^{n+1} = k \rho_{J}^{n+1} + (1 - k)\frac{\rho^*\lambda(\rho^*)}{\lambda^{n+1}} + kd^{n+1}, &  n \in \{0,1,\ldots\}.
\end{cases}
\end{equation}

We now define discrete version of ISS and ISS-Lyapunov function as follows 

\begin{definition}[Discrete ISS]\label{def:Disc_ISS}
	An equilibrium $ \rho^* \geq 0$ of the discrete closed-loop system \cref{eq:Disc_nonlocal_system1} is ISS in $L^2-$norm with respect to discrete disturbances  $ d^{n}$, $n \in \{1,2,\ldots\}$ if there exist positive real constants $\gamma_1 > 0$, $ \gamma_2 > 0 $ and $ \gamma_3 > 0 $ such that, for every initial condition $\rho_j^0$, $j \in \{1, \dots , J\}$, the solution $\rho_j^n $, $j \in \{1,\ldots,J\}$, $n \in \{0, 1, \ldots \} $ to the discrete closed-loop system \cref{eq:Disc_nonlocal_system1} satisfies
	\begin{equation}\label{eq:Disc_ISS_cond}
		\|\overrightarrow{\rho}^n \|_{L_{\Delta x}^2} \leq \gamma_2e^{-\gamma_1 t^n} \|\overrightarrow{\rho}^0 \|_{L_{\Delta x}^2} + \gamma_3 \max_{0 \leq s < n} \left(|d^s|\right), \; n \in \{1,2, \ldots\},
	\end{equation}
	where	
	\begin{equation*}
	\|\overrightarrow{\rho}^n \|^2_{L_{\Delta x}^2}:={\Delta x} \sum_{j=1}^{J} \left(\rho_j^n\right)^2,\quad n \in \{0, 1, \ldots \}. 
	\end{equation*} 
\end{definition}

\begin{definition}[Discrete ISS-Lyapunov function]\label{def:Disc_ISS_Lyapunov_fun}
	A discrete function $ \mathcal{L}^n > 0$, for all $n \in \{0,1, \ldots\} $ is said to be a discrete ISS-Lyapunov function for the discrete closed-loop system \cref{eq:Disc_nonlocal_system1} if 
	\begin{itemize}
		\item[(i)] there exist positive constants $ \alpha_1 > 0$ and $\alpha_2 > 0 $ such that for all $\rho_j^n $, $j \in \{1,\ldots,J\}$, $n \in \{0, 1, \ldots \} $
		\begin{equation}\label{eq:Disc_ISS_Lyfun1}
			\alpha_1\|\overrightarrow{\rho}^n \|_{L_{\Delta x}^2}^2 \leq \mathcal{L}^n \leq \alpha_2\|\overrightarrow{\rho}^n \|_{L_{\Delta x}^2}^2,
		\end{equation}
		\item[(ii)] there exist positive constants $ \eta > 0$ and $ \nu > 0 $ such that for all $\rho_j^n $, $j \in \{1,\ldots,J\}$, $n \in \{0, 1, \ldots \} $
		\begin{equation}\label{eq:Disc_ISS_Lyfun2}
			\frac{\mathcal{L}^{n+1} - \mathcal{L}^{n}}{\Delta t} \leq -\eta \mathcal{L}^{n} + \nu  \left(d^n\right)^2.
		\end{equation}
	\end{itemize}
\end{definition}

\begin{theorem}(Discrete ISS for $ \rho^* \geq 0 $)\label{thm:Disc_ISS_thm2}
	Assume that the CFL condition \cref{eq:CFL_cond} holds. For every $ \rho^* \geq 0 $, every $ k \in [0,1) $, every $ R > 0 $ and for every initial data $ \vec{\rho}^0 = (\rho_{0}^{0},\rho_{1}^{0}, \dots, \rho_{J}^{0})^\top $ satisfying  $ \rho_{j}^{0} \geq 0$, $j \in \{1, \dots, J\}$ and 
	\begin{equation}\label{eq:Disc_ISS_thm2_cond}
	\|\vec{\rho}^{0} - \rho^*\vec{e} \|_{L_{\Delta x}^2} \leq R, 
	\end{equation}
	where $ \vec{e} = \overbrace{(1,\ldots,1)^\top}^{J+1} $. Then, the solution $ \vec{\rho}^n = (\rho_{0}^{n},\rho_{1}^{n}, \dots, \rho_{J}^{n})^\top $ to the system \cref{eq:Disc_nonlocal_system1}  satisfies  $ \rho_{j}^{n} \geq 0$, $j \in \{0, \dots, J\}$, $n \in \{0,1,\ldots\}$ and the steady-state $ \rho^*$ of the discrete system \cref{eq:Disc_nonlocal_system1} is ISS in $L^2-$norm with respect to discrete disturbance function $ d^n$, $n \in \{1,2,\ldots\}$.
\end{theorem}
In order to analyze the ISS of the discrete system \cref{eq:Disc_nonlocal_system1} by the discrete Lyapunov method, we use the following transformation 
\begin{equation}\label{eq:Disc_transformation}
\tilde{\rho}_{j}^{n} = \rho_{j}^{n} - \rho^*, \; \widetilde{{W}}^n = {\Delta x}\sum_{j=1}^{J}\tilde{\rho}_{j}^{n},\quad \tilde{\lambda}_{\widetilde{{W}}}^n = \lambda\left(\rho^* + \widetilde{{W}}^n\right),\; \tilde{r}^n = \frac{\Delta t}{\Delta x} \tilde{\lambda}_{\widetilde{{W}}}^n,\; n \in \{0,1,\ldots\}.  
\end{equation}
For simplicity, we omit the symbol ``\textasciitilde" in \cref{eq:Disc_transformation} and discretize the system \cref{eq:Nonlocal_system3} as follows 
\begin{equation}\label{eq:Disc_nonlocal_system2}
\begin{cases}
\rho_{j}^{n+1} = \left( 1 - r^n\right) \rho_{j}^{n} + r^n \rho_{j-1}^{n},\; j \in \{1, \dots, J\},\; n \in \{0,1,\ldots\},\\
\rho_{0}^{n+1} = k \rho_{J}^{n+1} + (1 - k)\theta {W}^{n+1} + kd^{n+1}\;\text{with}\; \theta = \frac{\rho^*}{B + \rho^*},\;  n \in \{0,1,\ldots\},\\
r^n = \frac{\Delta t}{\Delta x}\lambda_{{W}}^n,\; n \in \{0,1,\ldots\},\\
\lambda_{{W}}^n = \lambda\left(\rho^* + {W}^n\right),\; n \in \{0,1,\ldots\},\\
{W}^n = {\Delta x}\sum_{j=1}^{J}{\rho}_{j}^{n} \geq - \rho^*,\; n \in \{0,1,\ldots\},\\
\lambda(s) = \frac{A}{B + s},\; s \geq 0. 
\end{cases}
\end{equation}
Thus, the assumption \cref{eq:Disc_ISS_thm2_cond} in \cref{thm:Disc_ISS_thm2} is now expressed as  
\begin{equation}\label{eq:Disc_ISS_thm2_cond2}
\|\vec{\rho}^{0} \|_{L_{\Delta x}^2} \leq R. 
\end{equation}  
\begin{proof}
	Note that the proof of \cref{thm:Disc_ISS_thm2} is an extension of the proof of Theorem 4.2 in \cite{Chen_2017}. Thus, some details of the proof can be found in \cite{Chen_2017}. Since the initial data $ \rho_{j}^{0} \geq 0$, $j \in \{0, \dots, J\} $, by the discrete system \cref{eq:Disc_nonlocal_system2} and the CFL condition \cref{eq:CFL_cond}, we have $ \rho_{j}^{n} \geq 0$, $j \in \{0, \dots, J\}$, $n \in \{0,1,\ldots\}$. Consider the following discrete Lyapunov function \cref{eq:Lyapunov_function} 
	\begin{equation}\label{eq:Disc_Lyapunov_function}
	\mathcal{L}^n = {\Delta x}\sum_{j = 1}^{J} (\rho_{j}^{n})^2 e^{-\beta x_j} + a(W^n)^2, \quad n \in \{0,1,\ldots\}, 	
	\end{equation}
	where $ \beta > 0 $ and $ 0 \geq a > -\frac{\beta}{e^{\beta}-1} $ are constants.	
	
	If we take 
	\begin{equation}\label{eq:Constant_a3} 
	a > - \left(\frac{{\Delta x}e^{\frac{\beta{\Delta x}}{2} }(e^{\beta} -1)}{e^{\beta{\Delta x}} -1}\right)^{-1},  
	\end{equation}
	then $ \mathcal{L}^n > 0 $ for all $ n \in \{0,1,\ldots\} $ and there exist positive constants $ C_1 > 0$ and $ C_2 > 0 $ such that  
	\begin{equation}\label{eq:Disc_Lyapunov_fun_approx1} 
	(W^n)^2 \leq C_1{\Delta x}\sum_{j = 1}^{J} (\rho_{j}^{n})^2 e^{-\beta x_j} \leq  \mathcal{L}^n \leq C_2 {\Delta x}\sum_{j = 1}^{J} (\rho_{j}^{n})^2 e^{-\beta x_j}, \quad n \in \{0,1,\ldots\}. 
	\end{equation}
	By using the discrete Lyapunov function \cref{eq:Disc_Lyapunov_function} and the discrete system \cref{eq:Disc_nonlocal_system2}, the time derivative of the Lyapunov function \cref{eq:Lyapunov_function} is approximated by  
	\begin{align}
	\frac{\mathcal{L}^{n+1} - \mathcal{L}^{n}}{{\Delta t}} &= \frac{{\Delta x}}{{\Delta t}} \sum_{j = 1}^{J} \left[\left(\rho_{j}^{n+1}\right)^2 - \left(\rho_{j}^{n}\right)^2 \right]e^{-\beta x_j}\notag\\
	& + \frac{a({\Delta x})^2}{\Delta t}\left[\left(\sum_{j = 1}^{J}\rho_{j}^{n+1}\right)^2 - \left(\sum_{j = 1}^{J}\rho_{j}^{n}\right)^2 \right],\notag\\
	&=  \left(e^{-\beta {\Delta x}} - 1\right){\lambda_{{W}}^n} \sum_{j = 1}^{J} \left(\rho_{j}^{n}\right)^2 e^{-\beta x_{j}} + e^{-\beta {\Delta x}} {\lambda_{{W}}^n}  \left[ \left(\rho_{0}^{n}\right)^2 - \left(\rho_{J}^{n}\right)^2e^{-\beta }\right]\notag\\	
	& + a {\Delta t}\left(\lambda_{{W}}^n\right)^2\left(\rho_{0}^{n} - \rho_{J}^{n}\right) + 2a\lambda_{{W}}^n \left(\rho_{0}^{n} - \rho_{J}^{n}\right)W^n,\; n \in \{0,1,\ldots\}.\label{eq:Disc_Lyapunov_fun_approx2} 
	\end{align}

	By using the boundary condition in \cref{eq:Disc_nonlocal_system2} and taking $ a \leq 0 $ as given by \cref{eq:Constant_a2} , we obtain that 
	\begin{align}
	\frac{\mathcal{L}^{n+1} - \mathcal{L}^{n}}{{\Delta t}} \leq&\; -\beta  e^{-\beta {\Delta x}}{\lambda_{{W}}^n} {\Delta x}\sum_{j = 1}^{J} \left(\rho_{j}^{n}\right)^2 e^{-\beta x_{j}}\notag\\
	& + e^{-\beta {\Delta x}}{\lambda_{{W}}^n} \left[ \left(k \rho_{J}^{n} + (1 - k)\theta {W}^{n} + kd^{n}\right)^2 - \left(\rho_{J}^{n}\right)^2e^{-\beta }\right]\notag\\	
	& + 2a\lambda_{{W}}^n \left(k \rho_{J}^{n} + (1 - k)\theta {W}^{n} + kd^{n} - \rho_{J}^{n}\right)W^n,\notag\\
	\leq&\; - \beta  e^{-\beta {\Delta x}}{\lambda_{{W}}^n} {\Delta x}\sum_{j = 1}^{J} \left(\rho_{j}^{n}\right)^2 e^{-\beta x_{j}} + k^2\left(1 + 2e^{-\beta {\Delta x}}\right) {\lambda_{{W}}^n}(d^{n})^2\notag\\
	& + {\lambda_{{W}}^n} b^n,\; n \in \{0,1,\ldots\},\label{eq:Disc_Lyapunov_fun_approx3} 
	\end{align}
	where 
	\begin{eqnarray*}
	b^n :=  e^{-\beta {\Delta x}} \left[ 2\left(k \rho_{J}^{n} + (1 - k)\theta {W}^{n}\right)^2  - \left(\rho_{J}^{n}\right)^2e^{-\beta }\right] - 2a(1 - k)\left(\rho_{J}^{n} - \theta {W}^{n}\right)W^n + a^2\left(W^n\right)^2.
	\end{eqnarray*}
	By substituting $ a $ and using convexity, the boundary term is simplified as follows 
	\small		
	\begin{align*}
	b^n = &\; e^{-\beta {\Delta x}}\left(2k^2 - e^{-\beta }\right)\left(\rho_{J}^{n}\right)^2 + 2(1 - k)e^{-\beta {\Delta x}}\theta^2 ({W}^{n})^2 \notag \\
	& - 2\theta(2k - e^{-\beta})\left(\rho_{J}^{n} - \theta {W}^{n}\right)W^n +  \left(\frac{2k - e^{-\beta}}{1 - k} \right)^2 \theta^2(W^n)^2,\notag\\
	=&\; \left(2k^2 - e^{-\beta }\right) e^{-\beta {\Delta x}}\left(\rho_{J}^{n} - \theta e^{-\beta {\Delta x}} {W}^{n} \right)^2\notag \\
	&+ \left[- \left(2k - e^{-\beta }\right) e^{\beta {\Delta x}} + 2(1 - k)e^{-\beta {\Delta x}} + 2(2k - e^{-\beta}) + \left(\frac{2k - e^{-\beta}}{1 - k} \right)^2\right] \theta^2(W^n)^2.\notag
	\end{align*}
	\normalsize
	For a given $ k \in [0,1) $ if one chooses $ \beta > 0 $ such that \cref{eq:Beta_condition1}  holds, then 
	\begin{equation}\label{eq:Disc_Lyapunov_fun_approx4} 
	b^n \leq \left[\left(2 - e^{\beta {\Delta x}}\right)\left(2k - e^{-\beta }\right)  + 2(1 - k)e^{-\beta {\Delta x}} + \left(\frac{2k - e^{-\beta}}{1 - k} \right)^2\right] \theta^2(W^n)^2.
	\end{equation}
	Using \cref{eq:Disc_Lyapunov_fun_approx1,eq:Disc_Lyapunov_fun_approx4}, we estimate \cref{eq:Disc_Lyapunov_fun_approx3}  as
	\begin{align}
	\frac{\mathcal{L}^{n+1} - \mathcal{L}^{n}}{{\Delta t}} &\leq  - f_{\Delta x}(\beta){\lambda_{{W}}^n} {\Delta x}\sum_{j = 1}^{J} \left(\rho_{j}^{n}\right)^2 e^{-\beta x_{j}}\notag\\
	&\quad + \frac{1}{2}\left(1 + 2e^{-\beta {\Delta x}}\right)e^{-\beta} {\lambda_{{W}}^n}(d^{n})^2,\; n \in \{0,1,\ldots\},\label{eq:Disc_Lyapunov_fun_approx5} 
	\end{align}
	where 
	\begin{align*}
	f_{\Delta x}(\beta) &:= \beta  e^{-\beta {\Delta x}}\\
	&\quad  - C_1\left[\left(2 - e^{\beta {\Delta x}}\right)\left(2k - e^{-\beta }\right)  + 2(1 - k)e^{-\beta {\Delta x}} + \left(\frac{2k - e^{-\beta}}{1 - k} \right)^2\right] \theta^2.
	\end{align*}
	Now, we take $ \beta = \beta(\rho^*,k) > 0 $ such that \cref{eq:Constant_a3,eq:Beta_condition1}  hold and 
	\begin{equation}\label{eq:Disc_Lyapunov_fun_approx6} 
	f_{\Delta x}(\beta) > 0.	
	\end{equation}
	Then, by using \cref{eq:Disc_Lyapunov_fun_approx1} , we obtain that  
	\begin{equation}\label{eq:Disc_Lyapunov_fun_approx7} 
	\frac{\mathcal{L}^{n+1} - \mathcal{L}^{n}}{{\Delta t}} \leq - \frac{f_{\Delta x}(\beta)}{C_1}{\lambda_{{W}}^n} \mathcal{L}^n + \frac{1}{2}\left(1 + 2e^{-\beta {\Delta x}}\right)e^{-\beta} {\lambda_{{W}}^n}\left(d^{n}\right)^2,\; n \in \{0,1,\ldots\}.
	\end{equation} 
	Recursively solving \cref{eq:Disc_Lyapunov_fun_approx7}, we obtain that 
	\begin{align}
	\mathcal{L}^{n+1} &\leq \left(1 - {\Delta t}\frac{f_{\Delta x}(\beta)}{C_1}{\lambda_{{W}}^n}\right) \mathcal{L}^n + \frac{1}{2}{\Delta t}\left(1 + 2e^{-\beta {\Delta x}}\right)e^{-\beta} {\lambda_{{W}}^n} \left(d^{n}\right)^2,\notag\\
	&\leq \prod_{m=0}^{n} \left(1 - {\Delta t}\frac{f_{\Delta x}(\beta)}{C_1}{\lambda_{{W}}^m}\right) \mathcal{L}^0\notag\\
	& + \frac{1}{2}\left(1 + 2e^{-\beta {\Delta x}}\right)e^{-\beta}{\Delta t}\sum_{m=0}^{n}{\lambda_{{W}}^m} \left(d^{m}\right)^2\prod_{r=m+1}^{n}\left(1 - {\Delta t}\frac{f_{\Delta x}(\beta)}{C_1}{\lambda_{{W}}^r}\right),\notag\\
	&\leq \exp\left(-\frac{f_{\Delta x}(\beta)}{C_1}{\Delta t}\sum_{m=0}^{n}\lambda_{{W}}^m \right)\mathcal{L}^0 \notag\\
	& + \frac{1}{2}\left(1 + 2e^{-\beta {\Delta x}}\right)e^{-\beta}{\Delta t}\sum_{m=0}^{n}{\lambda_{{W}}^m} \left(d^{m}\right)^2 \exp\left(-\frac{f_{\Delta x}(\beta)}{C_1}{\Delta t}\sum_{r=m+1}^{n}\lambda_{{W}}^r \right).\label{eq:Disc_Lyapunov_fun_approx8} 
	\end{align}
	By \cref{eq:Disc_nonlocal_system2,eq:Disc_ISS_thm2_cond2,eq:Disc_Lyapunov_fun_approx1,eq:Disc_Lyapunov_fun_approx8}, there exist positive constants $ C_3 > 0$ and $ C_4 > 0 $ such that for all $ n \in \{1,2,\ldots\} $,
	\begin{equation}\label{eq:Disc_Lyapunov_fun_approx9}
	-\rho^* \leq W^n \leq  \sqrt{\mathcal{L}^{n}} \leq \sqrt{C_3\mathcal{L}^{0} + C_4\max_{0 \leq s < n} \left(\left(d^s \right)^2\right)} \leq \sqrt{C_2C_3R^2 + C_4\max_{0 \leq s < n} \left(\left(d^s \right)^2\right)}.
	\end{equation}
	Let 
	\begin{equation}\label{eq:Disc_Lyapunov_fun_approx10}
	\sigma_2 := \min_{W^n} \lambda_{W}^n > 0,\quad \text{and} \quad \delta_2 := \max_{W^n} \lambda_{W}^n,\; n \in \{0,1,\ldots\}.  
	\end{equation}
	We use \cref{eq:Disc_Lyapunov_fun_approx10} in \cref{eq:Disc_Lyapunov_fun_approx8}  to obtain
	\begin{align}
	C_1\|\vec{\rho}^{n} \|_{L_{\Delta x}^2}^2 \leq \mathcal{L}^{n} &\leq e^{-\frac{\sigma_2f_{\Delta x}(\beta)}{C_1} t^{n}} \mathcal{L}^0 + \frac{1}{2}\left(1 + 2e^{-\beta {\Delta x}}\right)\frac{C_1\delta_2e^{-\beta}}{\sigma_2f_{\Delta x}(\beta)} \max_{0 \leq s < n} \left(\left(d^s \right)^2\right),\notag\\
	&\leq C_2e^{-\frac{\sigma_2f_{\Delta x}(\beta)}{C_1} t^{n}} \|\vec{\rho}^{0} \|_{L_{\Delta x}^2}^2 \notag\\
	&+ \frac{1}{2}\left(1 + 2e^{-\beta {\Delta x}}\right) \frac{C_1\delta_2e^{-\beta}}{\sigma_2f_{\Delta x}(\beta)} \max_{0 \leq s < n} \left(\left(d^s \right)^2\right),\; n \in \{1,2,\ldots\}.\label{eq:Disc_Lyapunov_fun_approx11}
	\end{align}	
	This concludes the proof of \cref{thm:Disc_ISS_thm2}. 
\end{proof}

\section{Numerical experiments}\label{sec:Section03}
In this section, we illustrate the theoretical results in \cref{sec:Section01,sec:Section02} by providing numerical computations of ISS of a scalar conservation law with nonlocal velocity and boundary measurement error. For this reason, we perform numerical computations of examples of a nonlocal conservation law with measurement error and compare results for ISS.     

We consider the closed-loop system \cref{eq:Nonlocal_system} with the given velocity function 
\begin{equation}\label{eq:Velocity_function_example}
\lambda(W(t)) = \frac{1}{1 + W(t)},\quad \text{with}\quad W(t) = \int_{0}^{1}\rho(t,x)dx,
\end{equation}
and rate of measurement error  
\begin{equation}\label{eq:measurementerrors}
d(t) = 2.4\times 10^{-3} \sin(t),\quad t \in (0, \infty). 
\end{equation} 

\subsection{Example 1}
\begin{figure}[h]
	\centering
	\includegraphics[width=\linewidth]{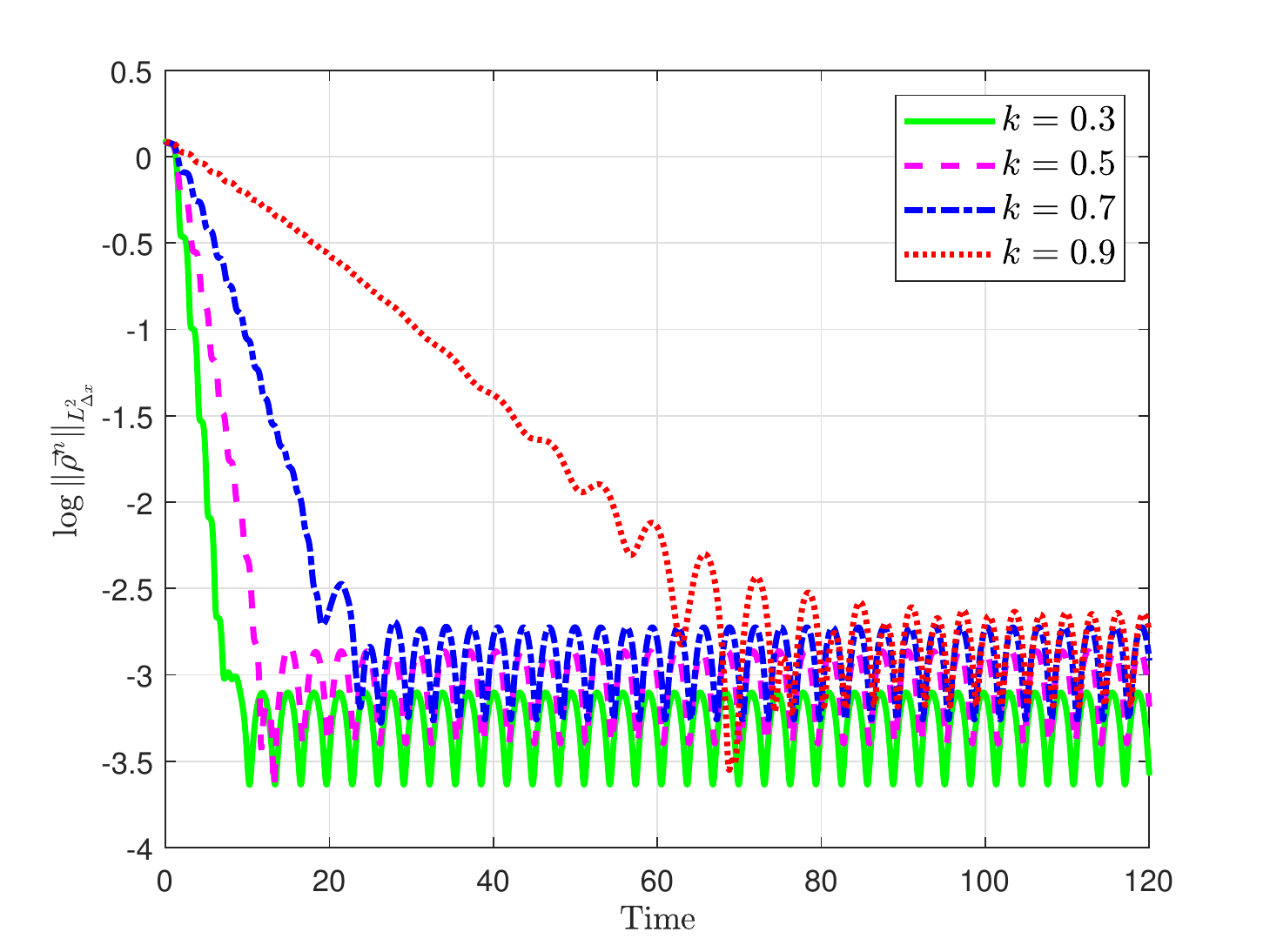}
	\caption{Comparison of Log-scale of $\|\vec{\rho}^n - \rho^*\vec{e}\|_{L_{\Delta x}^2}$ with CFL = 0.75  and $\rho^* = 0$.}
	\label{fig:nonlocal_trp_eqn_1a}
\end{figure}

Given an equilibrium solution $ \rho^* = 0 $, we set an initial condition $ \rho_0(x) = 1 + \sin(2\pi x)$ for $x \in [0,1]$. Moreover, for $ d \equiv 0 $ and in the sense of \cref{def:Disc_ISS_Lyapunov_fun} the decay rate of the Lyapunov function is obtained as follows 
\begin{equation}\label{eq:decay_rate_Lyapunov_fun}
	\eta:= \sigma_2 f_{\Delta x}(\beta),
\end{equation}
with $\displaystyle \sigma_2 := \min_{W(t)} \lambda(\rho^* + W(t)) $ and 
\begin{align*}
	f_{\Delta x}(\beta) &:= \beta  e^{-\beta {\Delta x}}\\
	&\quad  - \left[\left(2 - e^{\beta {\Delta x}}\right)\left(2k - e^{-\beta }\right)  + 2(1 - k)e^{-\beta {\Delta x}} + \left(\frac{2k - e^{-\beta}}{1 - k} \right)^2\right] \theta^2,
\end{align*}
where $ \theta = \frac{\rho^*}{1 + \rho^*} $, and $ k \in [0,1) $ and $ \beta > 0 $ are taken such that stability conditions hold. Besides, in the sense of \cref{def:Disc_ISS}, the discrete decay rates of the solution is given by $ \gamma_1 := 0.5\eta $. Then, we show the $ L^2-$error of the solution of the system \cref{eq:Nonlocal_system} and the discrete decay rates for two given CFL conditions 0.5 and 0.9 in \cref{tab:table_3}, respectively. Due to the artificial diffusion and the disturbance we observe only approximately the first--order of the numerical scheme. Furthermore, \cref{fig:nonlocal_trp_eqn_1a} shows the convergence of the solution of the system \cref{eq:Nonlocal_system} to the equilibrium for different values of $ k $. In \cref{fig:nonlocal_trp_eqn_1a}, we observe that as $ k $ increases the rate of decay of the Lyapunov function decreases due to the weaker control action. Furthermore, we observe that below the mesh accuracy of $\Delta x=10^{-3}$ no further decay is observed.

\begin{table}[h]
	\centering
	\begin{subtable}{0.48\textwidth}
		\centering
		\begin{tabular}{cccc}
			\hline
			$J$ & $ L^2- $error & order & $ \gamma_1 $\\
			\hline
			100   &      1.9171e-05   &      --    & 0.1270\\       
			200   &      1.1899e-05   &     0.6881 & 0.1274\\
			400   &      6.9631e-06   &     0.7730 & 0.1275\\
			800   &      3.7638e-06   &     0.8875 & 0.1276
			\\
			1600  &      1.5902e-06   &     1.2430 & 0.1277\\
			\hline
		\end{tabular}
		\caption{CFL = 0.5.}
		\label{tab:table_3a}
	\end{subtable}
	\hfill
	\begin{subtable}{0.48\textwidth}
		\centering
		\begin{tabular}{cccc}
			\hline
			$J$ & $ L^2- $error & order & $ \gamma_1 $\\
			\hline
			100   &      1.3831e-05   &      --    & 0.1270\\       
			200   &      8.1304e-06   &     0.7665 & 0.1274\\
			400   &      4.8604e-06   &     0.7423 & 0.1275\\
			800   &      2.8262e-06   &     0.7822 & 0.1276\\
			1600  &      1.1624e-06   &     1.2818 & 0.1277\\
			\hline
		\end{tabular}
		\caption{CFL = 0.9.}
		\label{tab:table_3b}
	\end{subtable}
	\caption{Comparison of $L^2- $error of the solution for number of grids $ J $ with $ \rho^* = 0 $, $ k = 0.3 $ and $ T = 10, $ where $ \gamma_1 $ is a rate of decay of the solution towards equilibrium.}
	\label{tab:table_3}
\end{table}

\subsection{Example 2}
\begin{figure}[h]
	\centering
	\includegraphics[width=\linewidth]{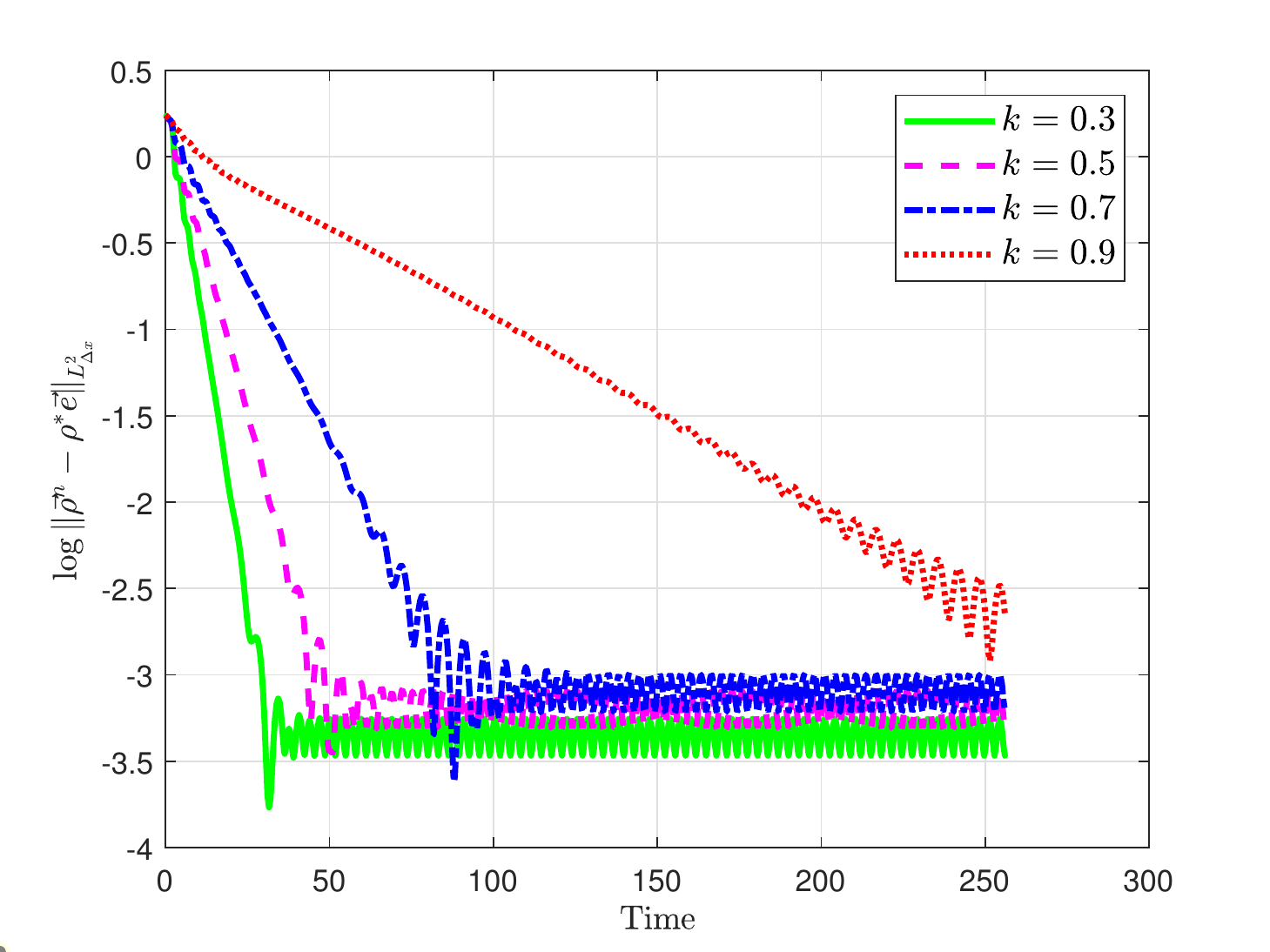}
	\caption{Comparison of Log-scale of $\|\vec{\rho}^n - \rho^*\vec{e}\|_{L_{\Delta x}^2}$ with CFL = 0.75 and $\rho^* = 1$.}
	\label{fig:nonlocal_trp_eqn_1b}
\end{figure}

We take an equilibrium solution $ \rho^* = 1 $ and an initial condition $ \rho_0(x) = 2 + 2\sin(2\pi x) \; x \in [0,1]$. We show similar results as above for the system \cref{eq:Nonlocal_system} with equilibrium $ \rho^* = 1 $ which are presented in \cref{tab:table_4,fig:nonlocal_trp_eqn_1b}.   Here, the first--order convergence of the scheme is visible. The observed decay rate $\gamma_1$ is smaller possibly due to the different equilibrium. 

\begin{table}[h]
	\centering
	\begin{subtable}{0.48\textwidth}
		\centering
		\begin{tabular}{cccc}
			\hline
			J & $ L^2- $error & order & $ \gamma_1 $\\
			\hline
			100   &      3.0916e-04   &      --    & 0.0244\\       
			200   &      1.5261e-04   &     1.0185 & 0.0244\\
			400   &      7.2438e-05   &     1.0750 & 0.0244\\
			800   &      3.1425e-05   &     1.2048 & 0.0244\\
			1600  &      1.0567e-05   &     1.5723 & 0.0244\\
			\hline
		\end{tabular}
		\caption{CFL = 0.5.}
		\label{tab:table_4a}
	\end{subtable}
	\hfill
	\begin{subtable}{0.48\textwidth}
		\centering
		\begin{tabular}{cccc}
			\hline
			J & $ L^2- $error & order & $ \gamma_1 $\\
			\hline
			100   &      2.8645e-04   &      --    & 0.0244\\       
			200   &      1.4299e-04   &     1.0024 & 0.0244\\
			400   &      6.9982e-05   &     1.0309 & 0.0244\\
			800   &      3.0215e-05   &     1.2117 & 0.0244\\
			1600  &      1.0128e-05   &     1.5769 & 0.0244\\
			\hline
		\end{tabular}
		\caption{CFL = 0.9.}
		\label{tab:table_4b}
	\end{subtable}
	\caption{Comparison of $L^2- $error of the solution for number of grids $ J $ with $ \rho^* = 1 $, $ k = 0.3 $ and $ T = 20, $ where $ \gamma_1 $ is a rate of decay of the solution towards equilibrium.}
	\label{tab:table_4}
\end{table}

\section{Conclusion}\label{Conclusion}

This paper considered Input-to-state stability (ISS) for a scalar conservation law with nonlocal velocity and boundary measurement error. A ISS-Lyapunov function is used to investigate conditions for ISS of an equilibrium for the scalar conservation law with nonlocal velocity and measurement error. Numerical study of a decay of ISS-Lyapunov function for such equations is analyzed. Finally, some examples are taken and numerical simulations are computed to illustrate the theoretical results. 

\section*{Acknowledgments} 
\addcontentsline{toc}{section}{Acknowledgments} 
The financial support of the DFG projects HE5386/18 and GO1920/10 is acknowledged.

\bibliographystyle{siamplain}
\bibliography{bibliography1} 

\end{document}